\documentclass[8pt]{article}
\usepackage{subcaption}
\usepackage{multicol}
\usepackage{multirow}
\usepackage{xcolor}
\usepackage{graphicx}
\usepackage{amssymb}
\usepackage{amsmath}

\usepackage{setspace}

\newtheorem{theorem}{Theorem}

\newtheorem{definition}{Definition}
\newtheorem{lemma}{Lemma}

\newcommand{\E}{{\mathsf E}}
\newcommand{\Pp}{{\mathsf P}}
\newcommand{\Var}{\mathsf{Var}}
\newcommand{\D}{\displaystyle}

\newcommand{\be}{\begin{equation}}
\newcommand{\ee}{\end{equation}}
\newcommand{\beq}{\begin{eqnarray}}
\newcommand{\eeq}{\end{eqnarray}}
\newcommand{\nbeq}{\begin{eqnarray*}}
\newcommand{\neeq}{\end{eqnarray*}}

\begin{document}


\title{Long-Term Behavior of Subordinated Branching Processes with Prevailing Emigration}
\author{George Yanev \\
 University of Texas Rio Grande Valley\\
Edinburg, Texas, U.S.A.\\
and \\
 Institute of Mathematics and Informatics\\
Sofia, Bulgaria\\
e-mail: {\tt george.yanev@utrgv.edu}}

\date{\today}

\maketitle

\begin{abstract}
This paper deals into the long-term behavior of subordinated critical branching processes with migration. We focus on scenarios where emigration is the dominant factor and introduce additional randomness in timing through a subordination mechanism, involving renewal processes. The key findings highlight how the initial population size and the interarrival mean time influence both asymptotic behavior of the non-extinction probability and corresponding Yaglom type limit theorems.  We also study an alternating regenerative process, when the population cycles between zero and positive states. This research complements previous studies for processes when immigration prevails over emigration.
\end{abstract}


subordination,  alternating regenerative process, Yaglom type limit 

theorems. 

\vspace{0.2 cm}{\bf MSC Classification:} 60J80, 60F05.

\section{Introduction}

Branching processes are fundamental stochastic models for describing the 
evolution of populations, epidemics, and other systems where entities 
reproduce and die over time. Classical models, such as the 
Bienaym\'e-Galton-Watson (BGW) process, assume discrete 
time and focus on reproduction as the sole mechanism 
driving population change. However, real-world populations are 
often subject to additional mechanisms such as migration 
(both immigration and emigration), external controls, and random 
environmental effects that influence their dynamics.

A particularly rich class of models arises when, in addition to the branching, 
migration is allowed. That is individuals
 or groups enter or leave the population at every generation.
  The interplay between reproduction and migration can lead to 
 a complex long-term behavior.
  In this paper, we focus on the regime 
  where emigration prevails over immigration, which leads to a net 
  outflow of individuals and a stronger tendency toward extinction.

Another important extension to the classical BGW process is the introduction of random time 
change, or subordination, when the process is observed at random 
times determined by the epochs of an independent 
renewal process. Generations are defined, but the process is updated only at random times according to a renewal process schedule. This approach captures additional sources of randomness reflecting variable environmental conditions, among others.
 Subordinated branching processes have found applications in fields ranging from biology to finance.

The main objective of this paper is to analyze the long-term 
behavior of critical BGW processes with migration, subordinated 
by a renewal process, in models with dominating emigration. We 
derive asymptotic results for the probability of non-extinction
 and Yaglom type conditional limit theorems. In our study,we highlight the influence of both the size of the initial population and the distribution of the interarrival times of the renewal process. We further extend the analysis to alternating regenerative processes, where the population alternates between positive and zero states, introducing additional complexity and leading to new limit theorems.

Our results extend previous studies (see \cite{P20}) on branching 
processes with migration and subordination, particularly in the context 
of heavy-tailed interarrival times and infinite mean number of progenitors.  The findings provide a comprehensive description of the asymptotic behavior of the processes when emigration is the dominant force. This offers new insights into the interplay between migration, random timing, and regenerative structure in branching systems. The proofs relay on previous results for processes with prevailing emigration and weighted renewal theory.

Commonly, a discrete-time stochastic process $\{Z_n\}_{n\in \cal I}$ has index set 
${\cal I}={\mathbb N_0}$.
However, $\cal I$ could be a set of integer numbers $j_1, j_2, j_3, \ldots $, realizations of a stochastic process $\{N(t)\}_{t\ge 0}$ with positive increments,
 so that $j_1 < j_2 < j_3 < ... $ Then a new
 process $\{Z_{N(t)}\}_{t\ge 0}$ can be constructed. This process is said to be subordinated to the underlying (base) process $\{Z_n\}_{n\in {\cal I}}$, while $\{N(t)\}$ is called time-change process.

The studies of subordinated branching processes have been advancing in two directions: considering different classes of branching processes and different types of subordination. In \cite{E96} a BGW process subordinated (randomly indexed) by a Poisson process was introduced  and studied as a model of daily stock prices. 
The more general model when the subordination is by a renewal process was studied in \cite{MMY09} and more recently in \cite {P20}, see also the references therein.


Let us have on a probability space $( \Omega ,{\cal F} , {\mathsf P})$ four
independent sets of non-negative, integer-valued, independent and identically distributed (i.i.d.) random
variables (r.v.s) representing the number of:  progenitors $Z_0$; offspring $\{ X_{n,i}\}_{n=
0}^\infty$ for $i\ge 1$; immigrants $\{ I_n\}_{n=0}^\infty$; and emigrants $\{ ( E_{fam,n}, E_{ind,n})\}_{n=0}^\infty$. Here $E_{fam,n}$ and $E_{ind,n}$ are the number of emigrating families and emigrating individuals, respectively.
Define  the migration component $\{ M_n\}_{n=0}^\infty$ of the process as follows:
\[ M_n = \left \{
\begin{array}{cl}
       {\D -\sum_{i=1}^{\D E_{fam,n}}X_{n,i}-  E_{ind,n}} & \mbox{with probability
\ {\it p},} \\
           0                    & \mbox{with probability \ {\it q},} \\ 
           
I_n       & \mbox{with probability \  {\it r}, \qquad  {\it p}+{\it q}+{\it
r}=1.}
            \end{array}
     \right. \]
     
\begin{definition} The BGW process with migration is defined by the recurrence
\[
Z_{n}=\left(\sum_{i=1}^{Z_{n-1}}X_{n,i}+M_{n}{\bf 1}
_{\{Z_{n-1}>0\}}\right)_+,\qquad n=1, 2, \ldots; \quad Z_{0}>0,
\]
where $(a)_+=\max\{0,a\}$ and ${\bf 1}_A$ is the indicator function.
\end{definition}\label{migration}

Any branching process with migration belongs to the general class of controlled ($\varphi$-branching) processes, see  \cite{GPY17} for discussion and references.
The process $\{ Z_n\}_{n\ge 0}$ is a homogeneous Markov chain, admitting
the following interpretation.  Upon the reproduction in the $n$th
generation, one of the three scenarios is possible:  {\bf (i)} the process evolves as an ordinary BGW process with probability $q$ ;
{\bf (ii)} ${\D
 E_{fam,n}}$ families and  ${ \D E_{ind,n}}$
individuals emigrate with probability $p$; {\bf (iii)}  $I_n$ immigrants join the population, given the generation is not empty, with probability $r$. Immigration is not permitted when the process is at zero, which makes this state absorbing. Clearly, Definition~1 includes as particular
cases processes with emigration (when $p=1$), the BGW process (when $q=1$), and processes with immigration outside zero (when $r=1$); see Figure 2. 

\begin{figure}[ht]
\includegraphics[height=2.6cm]{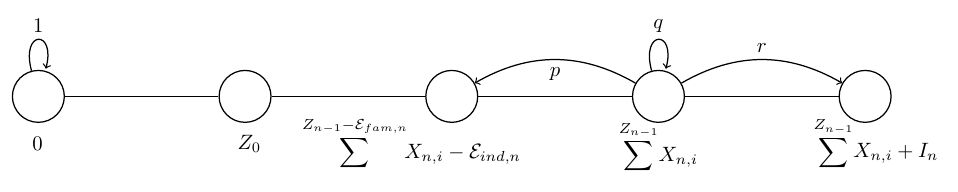}
\caption{Branching process with migration.}
\end{figure}

For notational convenience, we drop the subscripts $n$ and $i$ in $X_{n,i}$, $I_n$, $E_{fam,n}, E_{ind,n}$ and $M_n$ and use the notation $X, I, E_{fam}, E_{ind}$, and $M$, respectively.

Remarkably, the following single parameter $\theta$, involving both migration  and reproduction, plays a key role in the asymptotic behavior of  $Z_n$: 
\[
\theta:=\frac{2\E[M]}{\mathsf{Var}[X]}.
\]

We assume that $\E [M]<0$ yielding $\theta<0$. That is, we consider a 
BGW process with migration when emigration prevails over immigration. This complements the study in  \cite{P20} assuming $0<\theta<1$. We also need the following two sets of moments assumptions.

\vspace{0.3 cm}{\bf Assumptions A:} 
\begin{description}
  \item   -  (reproduction) $\E[X]=1$ and $0<\mathsf {Var} [X]<\infty$.
  \item   - (immigration) $\E[I]<\infty$.
 \item  - (emigration) $E_{fam}\le C_{fam}<\infty$ a.s. and 
$E_{ind}\le C_{ind}<\infty$ a.s.
\item  -  (migration) $\E[M]<0$. 
\end{description}

{\bf Assumptions B:}
 \[  
  \hspace {-2.7 cm}  
\begin{array}{ll}
    \E[X^2\log X]<\infty,
          & \mbox{if} \ -1<\theta<0.  \\
          & \\
 \E[X^2\log^2X]<\infty,   \quad  \E[I^2\log^2I]<\infty,	& \mbox{if } \qquad \quad \theta=-1. \\ 
    & \\
  \E[X^{1+|\theta|}]<\infty, \quad \quad \ \ \E [I^{1+|\theta|}]<\infty, & \mbox{if} \qquad \quad \ \theta<-1. 
   \end{array}
   \]

\section{Subordinated Process with Migration}

Let $\{N(t)\}_{t\ge 0}$ be a renewal (counting) process defined as follows

 \begin{definition} Let $ J_1, J_2 \ldots$ be a sequence of i.i.d. positive r.v.s (interarrival times).  Set $S_0=0$ and  $S_n=J_1+J_2+\ldots +J_n$. Then for $n\ge 0$
\[
N(t)=\max \{n: S_n\le t\}, \qquad t\ge 0.
\]
\end{definition}
The r.v. $N(t)$ counts the number of renewals (events) that have occurred by time $t$. We study  BGW processes with migration 
$\{Z_n\}_{n=0}^\infty$ subordinated by the renewal process 
$\{N(t)\}_{t\ge 0}$, which is independent of $\{Z_n\}_{n=0}^\infty$.

\begin{definition}
 A BGW process with migration $\{Z_n\}_{n=0}^\infty$ subordinated by the process $\{N(t)\}_{t\ge 0}$ is denoted by $\{Y(t)\}_{t\ge 0}$ and defined as follows: $Y(0)=Z_0$ and
$Y(t)=Z_{N(t)}$ for $t> 0$.
\end{definition}

The interpretation is: at time $S_0=0$, $Z_0$ progenitors (founders) produce offspring and accounting for migration, form the first generation of size $Z_1$ at time $S_1$. Next reproduction and migration occur at time $S_2$ and the size of the second generation is $Z_2$. The same procedure is repeated and the $n$th generation of size $Z_n$ is formed at time $S_n$. Clearly, if $t\in [ S_{n},S_{n+1})$, then $Y(t)=Z_{n}$. The population size remains constant between two consecutive renewal times. See Figure 1 for the particular case of subordinated BGW process.
Notice that $\{Y(t)\}_{t\ge 0}$ is not a Markov process unless $\{N(t)\}_{t\ge 0}$ is a homogeneous Poisson process.

Our objective is to investigate the limiting behavior of $Y(t)$ as $t\to \infty$. The process evolves until reaching the state zero (cemetery state). We are particularly interested in: (i) the rate of decay of the probability for non-extinction and (ii) the limiting distributions of $Y(t)$, appropriately normed and conditioned on $Y(t)>0$, i.e., given
non-extinction.

We are going to use the standard abbreviations ${\cal RV}_{-\alpha}$ and $\cal SV$ for the classes of regularly varying and slowly varying at infinity functions (cf. \cite{BGT86}).   We write $f\in {\cal{RV}}_{-\alpha}$ if $f(x)=x^{-\alpha}L(x)$ for $\alpha\in {\mathbb R}$, where $L(cx)/L(x)\to 1$ as $x\to \infty$ and $c\in {\mathbb{R}}$. The function $f$ is called regularly varying at infinity with $L$ being slowly varying at infinity function, $L\in \mathcal{SV}$ . 

The limiting behavior of $Y(t)$ as $t\to \infty$ is mainly regulated by two factors:  
 the average interarrival time $\mu$, say, and  the average number of progenitors $\E[Y(0)]$, each one being either finite or infinite. The various combinations of these give rise to a variety of limit theorems.
Consequently, we present the main results for $\{Y_t\}_{t\ge 0}$ in two groups. The first two theorems address the case $\mu<\infty$ and the  other two assume $\mu =\infty$ and $\Pp (J>x)= x^{-\rho}L_\rho(t)/\Gamma(1-\rho)$, $\rho\in (0,1)$. In the latter case, we deal with heavy-tailed distributed interarrival times.
 As far as the initial population size $Y(0)$ is concerned, we consider two possibilities: (i)  $\E[Y(0)]<\infty $ and (ii)  $\E[Y(0)]=\infty $ with
 $\Pp(Y(0)>x)\sim L_\gamma(x)x^{-\gamma}$ for
$\gamma\in (0,1)$ and large $x$.

\subsection{Finite Mean Interarrival Time: $\mu<\infty$}
 The results in this section show how the limiting properties of the discrete-time process $\{Z_n\}_{n=0}^\infty$ are transferred to the continuous-time process $\{ Y(t)\}_{t\ge 0}$ when $\mu<\infty$.
We begin with the probability of non-extinction. 

 \begin{theorem} \ \label{Thm1_Y}  
Suppose $\mu<\infty$ and Assumptions A hold. We have two claims:

$\rm (i)$ If $\E[Y(0)]<\infty$ and also Assumptions B hold, then
\[
{\mathsf P}(Y(t)>0)\sim L_\theta(t)
\left(\frac{t}{\mu} \right)^{-(1+|\theta|)}, \quad t\to \infty,
\]
where $L_\theta(t)\in \mathcal{SV}$ is given in Lemma~\ref{Zlemma1} in Section 4. 

$\rm (ii)$ If $\E[Y(0)]=\infty$ with $\Pp(Y(0)>x)\sim L_\gamma(x)x^{-\gamma}$, 
$\gamma\in (0,1)$,
then
\[
{\mathsf P}(Y(t)>0)\sim L_{\theta, \gamma}(t) \left(\frac{t}{\mu}\right)^{-\gamma},
\quad t\to \infty,
\]
where $L_{\theta, \gamma}(t)\in \mathcal{SV}$ is given in Lemma~\ref{Zlemma2} in Section 4.
\end{theorem}
Define the time to extinction (hitting zero) of  
$\{ Y(t) \}_{t\ge 0}$ by
\[
T= \inf\{ t \ge 1: Y(t)=0\} \quad \mbox{and} \quad  
\inf \emptyset =\infty .
\]
Clearly, Theorem~1 implies that: $\E [T]<\infty$ in (i) and $\E [T]=\infty$ in (ii) revealing the long-term effect of the heavy-tailed distribution of $\E[Z_0]$.

Next, we present an Yaglom type limit theorem. Let $2b:= \Var[X]$.
 
\begin{theorem}\ \label{Thm2_Y}
Suppose $\mu<\infty$ and Assumptions A hold.

$\rm (i)$ If $\E[Y(0)]<\infty$ and also Assumptions B hold, then
 \be \label{gamma1_finite_mu}
 \lim_{t\to \infty} {\mathsf P}\left(\frac{Y(t)}{b t/\mu}\le x | Y(t)>0\right)=
       {\mathsf P}(\xi \le x),
\ee
where $\xi$ is an exponential r.v., $\xi \sim {\rm Exp}(1)$.

$\rm (ii)$ If $\E[Y(0)]=\infty$ with ${\mathsf P}(Y(0)>0)\in {\mathcal{RV}}_{-\gamma}$,
$\gamma\in (0,1)$, then for $x\ge 0$
 \be \label{gamma2_finite_mu}
 \lim_{t\to \infty} {\mathsf P}\left(\frac{Y(t)}{b t/\mu}\le x | Y(t)>0\right)=
{\mathsf P}(\eta_{\theta, \gamma} \le x).
\ee
The limiting r.v. $\eta_{\theta,\gamma}$ has the following Laplace transform:
\be \label{LT}
\varphi(\lambda):= \E[e^{-\lambda \eta_{\theta,\gamma}}]=1+\theta\lambda^\gamma (1+\lambda)^{-\theta-\gamma}{\rm B}_{1/(1+\lambda)}(-\theta,1-\gamma),
\ee
where ${\rm B}_x(\cdot,\cdot)$ is the incomplete Beta function. 
\end{theorem}

The exponentially distributed limiting r.v. in (\ref{gamma1_finite_mu}) is the same as that appearing in the limit of the process with migration (see (\ref{lemma1ii}), Section 4). However, the normalization factor $bt/\mu$ in (\ref{gamma1_finite_mu}) is different due to the fact that, on average, the time is rescaled by $\mu$. It is worth mentioning that the same normalization appears in the critical Bellman-Harris process with mean life length $\mu$ (see \cite{AN72}, p.169). 

 Notice that 
$(1-\varphi(\lambda))/\lambda \sim c\lambda^{\gamma-1}$. Hence $\varphi(\lambda)\to 1$ as $\lambda\to 0+$, which makes it Laplace transform of a proper distribution. Since $\varphi(\lambda)\to \infty$ as $\lambda \to \infty$, this distribution does not have a mass at the origin.
It can also be shown (see \cite{YY97}) that in this case the tail behavior of the limiting distribution ($\sim x^{-\gamma}$) mimics that of $Y(0)$. Thus, we observe a long memory effect of the heavy-tailed $Z_0$ distribution. 

\subsection{ Infinite  Mean Interarrival Time: $\mu=\infty$}

Let us turn to the case $\mu=\E[J]=\infty$, assuming that the distribution of $J$ has a regularly varying heavy tail, i.e., for large $x$,
 \be \label{62}
\Pp (J>x)= x^{-\rho}L_\rho(t)/\Gamma(1-\rho), \qquad \rho\in (0,1).
\ee
One particular case of (\ref{62}) is the Fractional Poisson Process (cf. \cite{BGT86}). Its inter-arrival times' distribution has survival function
 $\mathsf{P}(J>t) \sim t^{-\rho}/(\lambda \Gamma(1-\rho))$, $\lambda>0$, 
(cf. \cite{Laskin2003}).

\begin{theorem}
 Suppose $\mu=\infty$ and (\ref{62}) is satisfied.
Let Assumptions A hold. Then:

$\rm (i)$ If $\E[Y(0)]<\infty$ and also Assumptions B hold, we have
\[
{\mathsf P}(Y(t)>0)\sim \E[T-1]L_\rho(t)t^{-\rho}, \quad t\to \infty.
\]
 $\rm (ii)$ If $\E[Y(0)]=\infty$ and $\Pp(Y(0)>x)\sim L_\gamma(x)x^{-\gamma}$, $\gamma\in (0,1)$,
then
\[
{\mathsf P}(Y(t)>0)\sim L_{\rho, \gamma}(t) t^{-\rho\gamma},\quad t\to \infty,
\]
where 
$
L_{\rho, \gamma}(t)=(1-\rho)^\gamma\Gamma^{1+\gamma}(1-\rho)\Gamma^{-1}(1-\rho\gamma)L_\rho(t)L_\gamma\left(\frac{t^\rho}{L_\rho(t)}\right) .
$
\end{theorem}

Under the assumptions of Theorem~3, $\E [T]=\infty$ and the non-extinction probability approaches zero at a slower rate than in Theorem~1 when $\mu<\infty$.

 To formulate the next Yaglom type result, we need the r.v. $\zeta_{\rho, \gamma}$ whose distribution function (d.f.) is
\be \label{zeta}
{\mathsf P}(\zeta_{\rho,\gamma} \le x) =  \frac{1}
{\E \left[\tau^{-\gamma}_\rho\right]}\int_0^x u^{-\gamma}\, dF_{\tau_\rho}(u)
     =  \frac{\Gamma(1-\rho\gamma)}{\Gamma(1-\gamma)}\int_0^x u^{-\gamma}\, dF_{\tau_\rho}(u). 
\ee
Here $\tau_\rho$ is a Mittag-Leffler r.v.
with Laplace transform 
$(1+\lambda^\rho)^{-1}$ and d.f. $F_{\tau_\rho}$ 
(see \cite{L98}). 
The r.v. $\zeta_{\rho,\gamma}$ is called a power-based version of $\tau_\rho$.

\begin{theorem}
 Suppose $\mu=\infty$ and (\ref{62}) is satisfied.
Let Assumptions A hold.

$\rm (i)$ If $\E[Y(0)]<\infty$ and also Assumptions B hold, then 
\be \label{thm2i} 
 \quad \lim_{t\to \infty} {\mathsf P}\left(\frac{Y(t)}{t^\rho/L_\rho(t)}\le x | Y(t)>0\right)= {\mathsf P}(\xi  \le x), \qquad x\ge 0,
\ee
where $\xi$ is an exponential r.v., $\xi \sim \rm{Exp}(1)$ and $L_\rho$ is given in (\ref{62}).

$\rm (ii)$ If $\E[Y(0)]=\infty$ and $\Pp (Y(0)>t)\in \cal{RV}_{-\gamma}$, $0<\gamma<1$,
then
\be \quad \lim_{t\to \infty} {\mathsf P}\left(\frac{Y(t)}{t^\rho/L_\rho(t)}\le x | Y(t)>0\right)= {\mathsf P}(\eta_{\theta, \gamma} \zeta_{\rho,\gamma} \le x) \qquad x\ge 0,
\ee
where the r.v.s $\eta_{\theta, \gamma}$ and $\zeta_{\rho, \gamma}$ are defined in (\ref{gamma2_finite_mu}) and (\ref{zeta}), respectively.
\end{theorem}

The limiting exponential law appears in both (\ref{gamma1_finite_mu}) and (\ref{thm2i}). However, the normalization factor $t^\rho/L_\rho(t)$ in (\ref{thm2i}) is sublinear. This means that the evolution of the process is slowed down compared to (\ref{gamma1_finite_mu}). Furthermore, it is known that in the infinite
  mean case, the renewal subordinator is asymptotically equivalent
   to a time-fractional subordinator (see \cite{Meerschaert2011}). This is consistent with the anomalous 
   scaling and Mittag-Leffler-type limit.

Table 1 summarizes our findings. Notice that
if $\mu<\infty$, then $Y(t)$ and $Z_n$ have the same limit, see Lemmas 2 and 3 for $Z_n$ in Section 4. If $\mu=\infty$ and the interarrival time has heavy-tailed distribution, then sublinear normalization is needed for $Y(t)$.

\begin{center}
\begingroup
\renewcommand{\arraystretch}{2} 
\begin{table}[ht]
\vspace{0.3cm} \begin{tabular}{|c|c|c|c|c|}
\hline
$\mu = \E[J]$ &  $\E[Y(0)]$  & ${\mathsf P}(Y(t)>0)$ & ${\mathsf P}\left(\frac{\D Y(t)}{\D A(t)}\le x | Y(t)>0\right)$ & Limit  \\ \hline
 \multirow{2}{*}{$\mu <\infty$} &  $\E[Y(0)]<\infty$    & ${\cal RV}_{-(1+|\theta|)}, \ \theta<0$ &  \multirow{2}{*}{$A(t)=bt/\mu$} & $\xi $  \\  
& $\E[Y(0)]=\infty$   & ${\cal RV}_{-\gamma}, \ \gamma\in (0,1)$ &  & $\eta_{\theta, \gamma}$   \\ \hline
     \multirow{2}{*}{$\mu =\infty$ } &   $\E[Y(0)]<\infty$  &  ${\cal RV}_{-\rho},\ \rho\in (0,1)$ & \multirow{2}{*}{$ A(t)= bt^\rho/L_\rho(t)$ } &  $\xi $\\ 
  &  $\E[Y(0)]=\infty$  &  ${\cal RV}_{-\rho \gamma}$ & & $\eta_{\theta, \gamma} \zeta_{\rho,\gamma}$   \\
    \hline
     \end{tabular}
     \caption{Summary of the results in Theorems 1-4.} 
     \end{table}
     \endgroup
\end{center}

\vspace{-1.5 cm}
\section{An Alternating Regenerative Process}

As Section 2 shows, the state zero is absorbing for the subordinated process with migration $\{Y(t)\}_{t\ge 0}$.
In this section we investigate a related process when zero is a reflective barrier. Namely, we construct and study a process, which spends some random sojourn time at zero before moving on to a positive state.

In general, the alternating regenerative processes are stochastic processes that 
alternate between two phases: ``down'' (failed or under repair) and ``up'' (operational), each lasting for a random duration. 
Regeneration points occur at the end of each cycle, 
where the process probabilistically restarts, making the future evolution independent of the past. 
Each cycle consists of a ``down'' period of random length $\tau^{(d)}$ and an ``up'' period of random length $\tau^{(u)}$,
so the cycle's length is $\tau =  \tau^{(d)} + \tau^{(u)}$.

Let $\{\tau^{(u)}_j : j = 1, 2, \ldots \}$ denote the sequence of up-periods during which the system (e.g., a machine) operates properly before a breakdown, and let 
 $\{\tau^{(d)}_j : j = 1, 2, \ldots \}$ denote the sequence of down-periods corresponding to time intervals during which the system is under repair or replacement. Assume that 
 these two sequences are mutually independent, and that the r.v.s within 
 each sequence are i.i.d. Define the $j$-th 
 cycle length as $\tau_j = \tau^{(d)}_j + \tau^{(u)}_j$ for $j = 1, 2, \ldots$. If the first up-period starts at time $\tau^{(d)}_1$, then the first breakdown occurs at time $\tau_1 = \tau^{(d)}_1 + \tau^{(u)}_1$. After the subsequent repair or replacement, 
  the system resumes operation at time $\tau_1 + \tau^{(d)}_2$, continuing until 
  the end of the second up-period at $\tau_2 = \tau_1 + \tau^{(d)}_2 + \tau^{(u)}_2$, 
  and so forth. Thus, each $\tau_j$ defines a renewal cycle. 

  Consider the renewal process $\{N_\tau(t)\}_{t\ge 0}$ generated by the sequence $\{\tau_j\}$ as follows
\[
N_\tau(t) = \max\{ n : S'_n \le t \}, \qquad \text{with} \qquad S'_n = \sum_{k=1}^n \tau_k, \quad S'_0 = 0.
\]
In addition, define the (spent time) process $\{\sigma(t)\}_{t \ge 0}$ by
\[
\sigma(t) := t - S'_{N_\tau(t)} - \tau^{(d)}_{N_\tau(t)+1}.
\]
The r.v. $\sigma^+(t) = \max\{ \sigma(t), 0 \}$ represents the  time  elapsed since the beginning of a up-period. Furthermore, let $\{Y_k(t)\}_{t\ge 0}$ for $k = 1, 2, \ldots$, be independent copies of the subordinated 
branching process with migration $\{Y(t)\}_{t\ge 0}$, 
with $\mathsf{P}(\tau_k \le t) = \mathsf{P}(T \le t)$, where $T$ is the extinction time of $\{Y(t)\}_{t\ge 0}$ defined in Section 2.1.

\begin{definition} \label{ARP} An alternating regenerative branching (ARB) process with migration $\{U(t): t\ge 0\}$ is defined by
\[
U(t)=
\left \{
\begin{array}{cl}
       Y_{N_\tau(t)+1}(\sigma(t)), & \mbox{if} \ \ \sigma(t)\ge 0 \\
       0 ,     & \mbox{if} \ \ \sigma(t)<0.
            \end{array}
     \right. 
\]
\end{definition}

Further on, since $\tau_k$ for $k=1,2,\ldots$ are i.i.d., we drop the subscript $k$ in $\tau_k$, $ \tau^{(d)}_k$, and $ \tau^{(u)}_k$. For the duration $ \tau^{(d)}$ of the sojourn period at zero assume that one of the following two conditions holds
\begin{subequations} \label{alpha}
\begin{equation}
\E \left[ \tau^{(d)}\right]<\infty
    \label{alpha1}
\end{equation}
\begin{equation} \label{alpha2}
 \E \left[ \tau^{(d)}\right]=\infty \quad  \mbox{and} \quad
    \Pp \left( \tau^{(d)} > t\right)\in {\cal{RV}}_{-\alpha}
    \quad 1/2<\alpha<1.
\end{equation}
\end{subequations}

Suppose $ \tau^{(d)}$ and $ \tau^{(u)}$ have non-lattice distributions, and $\Pp( \tau^{(d)}\le 0)=\Pp(\tau^   u\le 0)=0$. Consider the asymptotic ratio as $t\to \infty$ of the tails of the ``down'' (sojourn at zero) and ``up'' (branching) periods given by
\[
\Delta:=\lim_{t\to \infty}\frac{{\mathsf P}\left( \tau^{(d)}>t\right)}{{\mathsf P}\left(  \tau^{(u)}>t\right)}=\lim_{t\to \infty}\frac{{\mathsf P}\left( \tau^{(d)}>t\right)}{{\mathsf P}\left( Y(t)>0) \right)}, \qquad 0\le \Delta\le \infty.
\]

\noindent In our study we deal with both $\E[ \tau^{(d)}]<\infty$ and $\E[ \tau^{(d)}]=\infty$
Note, if the distribution of the sojourn period has heavier tail that, the distribution of the branching period, then $\Delta=0$.
Figure 2 illustrates the evolution of the ARB process with migration.
\begin{figure}[ht]
\includegraphics[height=5cm, width=13cm]{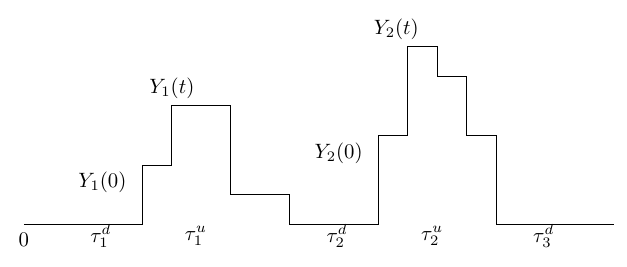}
\caption{Alternating branching process with migration.}
\end{figure} 

\subsection{Finite Mean Interarrival Time: $\mu<\infty$}

Here we assume that both $\mu<\infty$ and $\E Y(0)=\infty$ with
\be \label{infiniteY0}
\Pp(Y(0)>x)\sim L_\gamma(x)x^{-\gamma}, \quad 1/2<\gamma<1.
\ee
Note that if (\ref{alpha2}) holds, then (\ref{infiniteY0}) and Theorem~1(ii) yield
\[
\Delta:=\lim_{t\to \infty}\frac{L_\alpha(t)t^{-\alpha}}{ L_{\theta, \gamma}(t)(t/\mu)^{-\gamma} }, \qquad 1/2<\alpha, \gamma <1.
\]
 
\begin{theorem}\label{arp_thm1} 
    Suppose $\mu<\infty$ and $\E[Y(0)]=\infty$ with (\ref{infiniteY0}). Let Assumptions A hold.
    
$\rm (i)$ If, additionally, one of (\ref{alpha1}) or (\ref{alpha2}) is true and $0\le \Delta<\infty$, then 
\[
\lim_{t\to \infty}\Pp\left(\frac{U(t)}{bt/\mu}\le x\right) 
  = 
    \left\{
    \begin{array}{ll}
\frac{\D 1}{\D \Delta+1} (\Delta + \Pp(\eta_{\theta, \gamma} \beta_{\gamma, 1-\gamma} \le x)), &  \mbox{if} \ x>0, \\
 &  \\
\frac{\D \Delta}{\D \Delta+1}, & \mbox {if} \ \  x=0,
    \end{array}
        \right. 
\]
where 
$\eta_{\theta, \gamma}$ is the limiting r.v. in Theorem 2(ii), independent 
of  $\beta_{\gamma, 1-\gamma}\sim {\rm Beta}(\gamma, 1-\gamma)$. 

$\rm (ii)$ If, additionally, (\ref{alpha2}) is true and  $\Delta=\infty$, then for $x\ge 0$
\[
\lim_{t\to \infty}\Pp\left(\frac{U(t)}{bt/\mu}\le x\ \Big|\ U(t)>0\right)  
 =  
\Pp(\eta_{\theta, \gamma} \beta_{\alpha, 1-\gamma}\le x),
\] 
where 
$\eta_{\theta, \gamma}$ is the limiting r.v. in Theorem 2(ii),  independent of  $\beta_{\alpha, 1-\gamma}\sim {\rm Beta} (\alpha, 1-\gamma)$. 
\end{theorem}

\subsection{Infinite  Mean Interarrival Time: $\mu=\infty$}
In the next two theorems we assume $\mu=\infty$ with $\Pp (J>t)\in \cal{RV}_{-\rho}$ for $0<\rho<1$ 
and $\E [Y(0)]<\infty$ (Theorem 6) or $\E [Y(0)]=\infty$ (Theorem 7) with the additional condition $\Pp (Y(0)>x)\in \cal{RV}_{-\gamma}$ for $0<\gamma<1$.
Note that if (\ref{alpha2}) holds, then (\ref{infiniteY0}) and Theorem~3(i) yield
\[
\Delta:=\lim_{t\to \infty}\frac{L_\alpha(t)t^{-\alpha}}{ E[T-1]L_{\rho}(t)t^{-\rho} }, \qquad 1/2<\alpha, \rho <1.
\]
    
\begin{theorem} \label{arp_thm2}
    Suppose $\mu=\infty$ with
    ${\mathsf P}(J>t)\in {\cal{RV}}_{-\rho}$, $1/2<\rho<1$ and $\E[Y(0)]<\infty$. Let both Assumptions A and B hold . 

$\rm (i)$ If, additionally, one of (\ref{alpha1}) or (\ref{alpha2}) is true and $0\le \Delta<\infty$, then 
\[
\lim_{t\to \infty}\Pp\left(\frac{U(t)}{t^\rho/L_\rho(t)}\le x\right) 
  = 
    \left\{
    \begin{array}{ll}
\frac{\D 1}{\D \Delta+1} (\Delta + \Pp(\xi \beta^\rho_{\rho, 1-\rho} \le x))  &  \mbox{if} \ x>0, \\
 &  \\
\frac{\D \Delta}{\D \Delta+1} & \mbox {if} \ \  x=0,
    \end{array}
        \right. 
\]
where 
$\xi\sim {\rm Exp}(1)$, independent of  $\beta_{\rho, 1-\rho}\sim {\rm Beta} 
(\rho, 1-\rho)$. The random variable $\xi\beta^\rho_{\rho, 1-\rho}/(\Delta+1)$ has mean 
\[
\frac{1}{\Delta+1}\E [ \xi \beta^\rho_{\rho, 1-\rho} ]=\frac{\Gamma(2\rho)}{(\Delta+1)\Gamma^2(\rho)}.
\]

$\rm (ii)$ If, additionally, (\ref{alpha2}) is true and  $\Delta=\infty$, then for $x\ge 0$
\[
\lim_{t\to \infty}\Pp\left(\frac{U(t)}{t^\rho/L_\rho(t)}\le x\ \Big|\ U(t)>0\right)  
 =  \Pp(\xi \beta^\rho_{\alpha, 1-\rho} \le x),
\] 
where $\xi\sim {\rm Exp}(1)$, independent of  $\beta_{\alpha, 1-\rho}\sim {\rm Beta} 
(\alpha, 1-\rho)$. The limit $\xi \beta^\rho_{\alpha, 1-\rho}$  has mean 
\[
\E [\xi \beta^\rho_{\alpha, 1-\rho} ]=
\frac{\D {\rm B}(\alpha+1-\rho, \rho)}{\D {\rm B}(\alpha, \rho)}.
\]
\end{theorem}

Note that if in (\ref{alpha2}) is true, then (\ref{infiniteY0}) and Theorem~3(ii) yield
\[
\Delta:=\lim_{t\to \infty}\frac{L_\alpha(t)t^{-\alpha}}{ L_{\rho, \gamma}(t)t^{-\rho\gamma} }, \quad 1/2<\alpha, \rho \gamma <1, \quad 0<\rho, \gamma<1.
\]

\begin{theorem} \label{arp_thm3} 
    Suppose both $\mu$ and $\E Y(0)$ are infinite with $1/2<\rho \gamma<1$ and
    \[
    {\mathsf P}(J>t)\in {\cal{RV}_{-\rho}},\quad  0<\rho<1 \quad 
\mbox{and} \quad \Pp(Y(0)>x) \in {\cal{RV}}_{-\gamma},\quad 0<\gamma<1.
\]
Let also Assumptions A hold.

$\rm (i)$ If, additionally, one of (\ref{alpha1}) or (\ref{alpha2}) is true and $0\le \Delta<\infty$, then 
\[
\lim_{t\to \infty}\Pp\left(\frac{U(t)}{t^\rho/L_\rho(t)}\le x\right) 
  = 
    \left\{
    \begin{array}{ll}
\frac{\D 1}{\D \Delta+1} (\Delta + \Pp(\eta_{\theta, \gamma} \zeta_{\rho, \gamma} \beta^\rho_{\rho\gamma, 1-\rho\gamma}\le x) &  \mbox{if} \ x>0, \\
 &  \\
\frac{\D \Delta}{\D \Delta + 1} & \mbox {if} \ \  x=0.
    \end{array}
        \right. 
\]
Here $\eta_{\theta, \gamma}$ is the limiting r.v. defined in Theorem 2(ii) and $\zeta_{\rho, \gamma}$ is 
the power-based Mittag-Leffler r.v. defined in Theorem 3(ii). All random variables $\eta_{\theta, \gamma}$, $\zeta_{\rho, \gamma}$, and $\beta_{\rho\gamma, 1-\rho\gamma}\sim {\rm Beta}(\rho\gamma, 1-\rho\gamma)$  
are independent. 

$\rm (ii)$ If, additionally, (\ref{alpha2}) is true and  $\Delta=\infty$, then for $x\ge 0$
\[
\lim_{t\to \infty}\Pp\left(\frac{U(t)}{t^\rho/L_\rho(t)}\le x\ \Big|\ U(t)>0\right)  
 =  \Pp(\eta_{\theta, \gamma} \zeta_{\rho, \gamma} \beta^\rho_{\alpha, 1-\rho\gamma}\le x),
\] 
Here $\eta_{\theta, \gamma}$ is the limiting r.v. defined in Theorem 2(ii) and $\zeta_{\rho, \gamma}$ is 
the power-based Mittag-Leffler r.v. from Theorem 3(ii).  All r.v.s$\beta_{\alpha, 1-\rho\gamma}\sim {\rm Beta}(\alpha, 1-\rho\gamma)$, $\eta_{\theta, \gamma}$,and $\zeta_{\rho, \gamma}$ 
 are independent. 
\end{theorem}

The Beta distribution appears in the above theorems because, in the limit, the proportion of time spent in the current renewal cycle (relative to the total length of the renewal cycle) converges to a Beta law. Specifically, for $\rho\in(0,1)$, the limiting distribution of the normalized residual life is Beta$(\rho,1-\rho)$. In our context, the r.v. $\xi\beta_{\rho,1-\rho}^\rho$ in the limit, describes the scaled contribution of the ``current'' renewal cycle of the process $U(t)$ as $t\to\infty$. Here $\xi\sim {\rm Exp}(1)$ is an independent r.v., and $\beta_{\rho,1-\rho}$ captures the random proportion of the cycle. 

In Table~2 we collect the results for ARB process with migration under various regimes of the model's parameters. Additionally, Figure 3 and Figure 4 provide visual summaries setting $\theta=-0.2$, $\gamma=\rho=0.9$, and $\alpha=0.6$.

\begin{center}
\begingroup
\renewcommand{\arraystretch}{2} 
\begin{table}[ht]
\begin{tabular}{|c|c|c|c|c|}
\hline
$\mu = \E[J]$ &  $\E[Y(0)]$  & $\Delta$ & Probab.& Limit  \\ \hline
 \multirow{2}{*}{$\mu <\infty$} &  \multirow{2}{*}{$ \frac{1}{2}<\gamma<1$ }   & $[0,\infty)$&  $\Pp\left(\frac{U(t)}{bt/\mu}\le x\right)$ & $\eta_{\theta,\gamma} \beta_{\gamma, 1-\gamma} $  \\  
&   & $\frac{1}{2}<\alpha<1$ & $\Pp\left(\frac{U(t)}{bt/\mu}\le x | U(t)>0\right)$ & $\eta_{\theta,\gamma} \beta_{\alpha, 1-\gamma}$   \\ \hline
 $\mu =\infty$ &  \multirow{2}{*}{$ (0,\infty)$ }   & $[0,\infty)$&  $\Pp\left(\frac{U(t)}{t^\rho/L_\rho(t)}\le x\right)$ & $\xi \beta^\rho_{\rho, 1-\rho}$  \\  
$\frac{1}{2}<\rho<1$ &   & $\frac{1}{2}<\alpha<1$ & $\Pp\left(\frac{U(t)}{t^\rho/L_\rho(t)}\le x | U(t)>0\right)$ & $\xi \beta^\rho_{\alpha, 1-\rho} $   \\ \hline
 $\mu =\infty$ &  \multirow{2}{*}{$ \frac{1}{2}<\rho \gamma<1$ }   & $[0,\infty)$&  $\Pp\left(\frac{U(t)}{t^\rho/L_\rho(t)}\le x\right)$ & $\eta_{\theta, \gamma} \zeta_{\rho, \gamma} \beta^\rho_{\rho\gamma, 1-\rho\gamma} $   \\  
$0<\rho<1$ &  $0<\gamma<1$ & $\frac{1}{2}<\alpha<1$ & $\Pp\left(\frac{U(t)}{t^\rho/L_\rho(t)}\le x | U(t)>0\right)$ & $\eta_{\theta,\gamma} \zeta_{\rho, \gamma} \beta^\rho_{\alpha, 1-\rho\gamma}$   \\ \hline
     \end{tabular}
     \caption{Summary of the results in Theorems 5-7.} 
     \end{table}
     \endgroup
\end{center}

\begin{figure}[h]
\centering
\begin{minipage}{0.5\textwidth}
  \hspace{-1.5 cm}\includegraphics[height=14cm, width=9 cm]{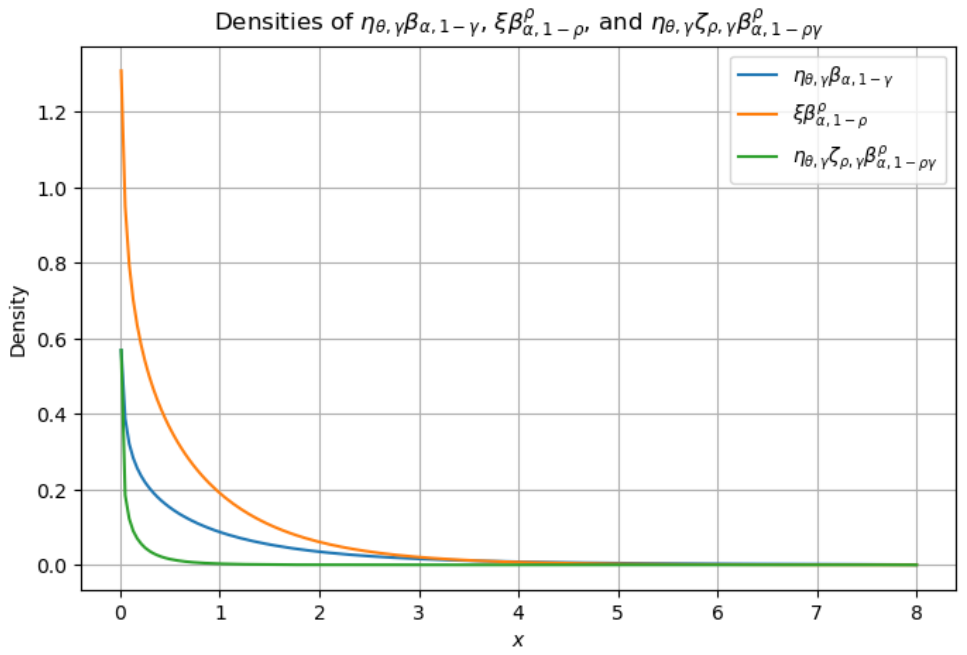}
 \vspace{-7.7 cm} \caption{Densities of ARB process \\ limits when $\Delta=\infty$.}
\end{minipage}\hfill
\begin{minipage}{0.5\textwidth}
 \hspace{-1 cm} \includegraphics[height=14cm, width=9 cm]{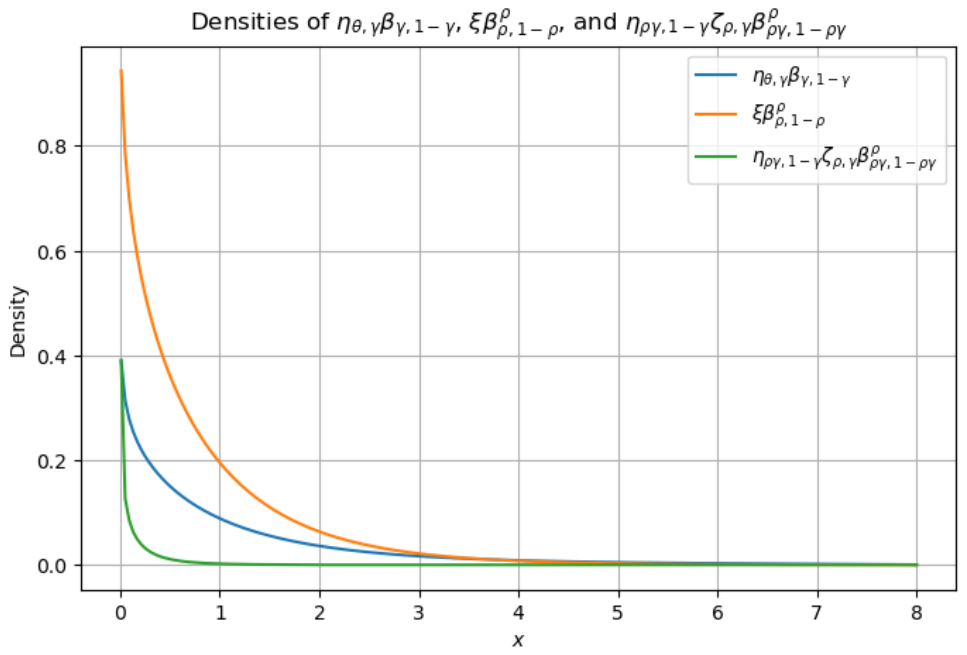}
 \vspace{-7.7 cm} \caption{Densities of the ARB process limits when $\Delta=0$.}
\end{minipage}
\end{figure}

\section{Auxiliaries}
For reader's convenience we introduce here some notations and formulate a few lemmas, which will be used in the proofs of the theorems in the next section.

\subsection{Processes with Migration}

In this subsection, we present some known results for the process with migration $\{Z_n\}_{n=0}^\infty$ when $\theta<0$. As before we use the notation $2b={\Var}[X]$.

\begin{lemma}[\cite{YY04}] \label{Zlemma1}
Let Assumptions A and B hold. If $\E[Z_0]<\infty$,
then 
\[
{\mathsf P}(Z_n>0)\sim L_\theta(n)n^{-(1+|\theta|)}, \qquad n\to \infty,
\]
where $L_\theta\in \cal{SV}$. Furthermore, with $\xi \sim {\rm Exp}(1)$
 \be \label{lemma1ii}
 \lim_{n\to \infty} {\mathsf P}\left(\frac{Z_n}{b n}\le x | Z_n>0\right)= {\mathsf P}(\xi \le x),\qquad x\ge 0.
\ee
\end{lemma}

In Lemma \ref{Zlemma1} we have ${\mathsf P}(Z_n>0) \sim L_\theta(n)n^{-(1+|\theta|)}$, which is a faster rate of decay than in the BGW process due to the dominating emigration. However, (see \cite{YY04})
\be \label{bn}
\E[Z_n\ |\ Z_n>0]=\frac{\E[Z_n]}{{\mathsf P}(Z_n>0)} \sim \frac{bL_\theta(n)n^{-|\theta|}}{L_\theta(n)n^{-(1+|\theta|)}}\sim bn.
\ee
This explains the normalization in (\ref{lemma1ii}). It turns out that, despite the different decay rates of their non-extinction probabilities, both $\{ Z_n\}_{n=0}^\infty$ and BGW processes, under the same linear normalization, share the same exponential limiting distribution on their non-extinction paths.

 Let, instead of $\E[Z_0]<\infty$, assume $\E[Z_0]=\infty$ and $\Pp(Y(0)>x)\sim L_\gamma(x)x^{-\gamma}$, 
$\gamma\in (0,1)$.
 Under these assumptions, the next lemma shows that, the no-extinction probability and the tail of of the distribution of $Z_0$ approach zero at the same rate.

\begin{lemma}[\cite{YY97}] \label{Zlemma2}
Suppose  $\E[Z_0]=\infty$ and
$\Pp(Y(0)>x)\sim L_\gamma(x)x^{-\gamma}$ where
$\gamma\in (0,1)$.  If Assumptions A hold, then
\be \label{lemma2i}
{\mathsf P}(Z_n>0) \sim L_{\theta,\gamma}(n)n^{-\gamma} \quad \mbox{where}\quad L_{\theta,\gamma}(n) =\frac{b^{-\gamma}(1-\theta-\gamma)^{-1}}{{\rm B}(1-\theta, 1-\gamma)}L_\gamma(n)
\ee
and ${\rm B}(\cdot, \cdot)$ is the Beta function.
Furthermore,
 \be \label{lemma2ii}
 \lim_{n\to \infty} {\mathsf P}\left(\frac{Z_n}{b n}\le x | Z_n>0\right)= {\mathsf P}(\eta_{\theta, \gamma} \le x),\qquad x\ge 0.
\ee
The Laplace transform of the limit r.v. $\eta_{\theta, \gamma}$ is given in (\ref{LT}).
\end{lemma}


In Lemma \ref{Zlemma2} the rate of decay of the survival probability (\ref{lemma2i}) depends on the parameter $\gamma$. That is, we again observe the long-term effect of the the assumption $\E [Z_0]=\infty$, which changes the limit r.v. However, in both (\ref{bn}) and (\ref{lemma2ii}) the normalization is $bn$.

\subsection{Weighted Renewal Functions}
In this subsection we collect known results from the weighted renewal processes \cite{MO06}. 

\begin{definition} Let $\{w_k\}_{k=0}^\infty$  be a sequence of non-negative constants (weights). A weighted renewal function $H$ is defined as 
\[
H(t)=\sum_{k=0}^\infty w_k {\mathsf P}(N(t)\ge k), \qquad  t\ge 0.
\]
\end{definition}
Denote
$w_0={\mathsf P}(Z_0=0)$, $w_k={\mathsf P}(Z_k=0)-{\mathsf P}(Z_{k-1}=0), \quad k=1, 2, \ldots $
By the law of total probability and the independence of $N(t)$ and $Z_k$, we have 
\beq \label{weights}
{\mathsf P}(Y(t)=0) & = & \sum_{k=1}^\infty {\mathsf P}(N(t)=k){\mathsf P}(Z_k=0) \\
        & = & \sum_{k=1}^\infty {\mathsf P}(N(t)=k)\sum_{n=1}^k w_n \nonumber \\
        & = & \sum_{k=1}^\infty w_k {\mathsf P}(N(t)\ge k). \nonumber 
\eeq
Denote $W(x)=\sum_{k=0}^{\lfloor{x}\rfloor} w_k={\mathsf P}(Z_{\lfloor{x}\rfloor} =0)$.

\begin{lemma}[Theorem 46(i), \cite{MO06}] \label{l1} Let $\mu<\infty$ and $L_\gamma(t)\in \cal{SV}$ 
 If $w_{k}\geq 0$ and $W(t)\uparrow 1$ as $t\to \infty$, then  for $0<\gamma<1$
\[
1-W(t)\sim L_\gamma(t)t^{-\gamma}\quad \mbox{ if and only if }\quad  1-H(t)\sim L_\gamma(t)\left(\frac{t}{\mu}\right)^{-\gamma}, \quad t\to \infty.
\]
\end{lemma}
Note also that (\cite{MMY10}) if $\mu<\infty$, then for $x\ge 1$
\be  \label{Dlimit}
\lim_{t\to\infty}{\mathsf P}\left( \frac{N(t)}{t/\mu}\le x | Z_{N(t)}>0\right)= \mathbf{1}_{\{x\ge 1 \}}.
\ee
We turn to the case $\mu=\infty$ assuming for $0<\rho<1$
\[
{\mathsf P}(J>t) \sim L_{\rho}(t)t^{-\rho},\qquad  t\to \infty,
\]
where $L_\rho\in {\cal SV}$. Recall that
$
T= \inf\{ t \ge 1: Y(t)=0\} \quad \mbox{and} \quad  
\inf \emptyset =\infty .
$
Since ${\mathsf P}(T=k)={\mathsf P}(Z_k=0)-{\mathsf P}(Z_{k-1}=0)$, we can rewrite (\ref{weights}) as
\be \label{Y and H}
{\mathsf P}(Y(t)=0) = \sum_{k=0}^\infty {\mathsf P}(T=k){\mathsf P}(N(t)\ge k)= H(t).
\ee
The next lemma summarizes some asymptotic properties of ${\mathsf P}(Y(t)>0)=1-H(t)$. 

\begin{lemma}[Theorem 1, \cite{OMM09}] \label{Thm1_12}
Suppose that 
\[
{\mathsf P}(J>x)\in {\cal{RV}}_{-\rho}, \quad 0<\rho <1 \quad \mbox{and}\quad {\mathsf P}(T>t)\in {\cal{RV}}_{-\gamma}, \quad \gamma \ge 0.
\]
$\rm (i)$ If $\gamma \ge 1$ and $\E[T]<\infty$, then
\[
\lim_{t\to \infty}\frac{1-H(t)}{{\mathsf P}(J>t)}=\E[T-1].
\]
$\rm (ii)$ If $0\le \gamma <1$ and $\E[T]=\infty$, then
\[
\lim_{t \to \infty}\frac{1-H(t)}{\D {\mathsf P}(T>t^\rho / L_\rho(t))}=\frac{\Gamma^\gamma(2-\rho)\Gamma(1-\rho)}{\Gamma(1-\rho\gamma)}.
\]
\end{lemma}

\noindent The last lemma establishes the limit of $N(t)$ given non-extinction of the $\{Y(t)\}_{t\ge 0}$.

\begin{lemma}[Theorem 2, \cite{OMM09}] \label{Thm2_12} 
Suppose that 
\[
{\mathsf P}(J>x)\in {\cal{RV}}_{-\rho}, \quad 0<\rho <1 \quad \mbox{and}\quad {\mathsf P}(T>t)\in {\cal{RV}}_{-\gamma}, \quad \gamma \ge 0.
\]
$\rm (i)$ If $\gamma \ge 1$ and $\E[T]<\infty$, 
then
\[
\lim_{t\to \infty}
{\mathsf P}\left(\frac{N(t)}{t^\rho / L_\rho(t)}\le  x\ |\ T> N(t)\right) = \mathbf{1}_{\{x\ge 1 \}}.
\]
$\rm (ii)$ If $0\le \gamma<1$ and $\E[T]= \infty$, then, with r.v. $\zeta_{\rho, \gamma}$ as given in 
(\ref{zeta}),
\[
\lim_{t\to \infty}{\mathsf P}\left(\frac{N(t)}{t^\rho / L_\rho(t)}\le x\ |\ T> N(t)\right)
= {\mathsf P}(\zeta_{\rho, \gamma}\le x).
\]
\end{lemma}
 
\section{Proofs of the Theorems for $\{Y(t)\}_{t\ge 0}$}

\begin{lemma} \label{integral}
    The following relation holds true as $x>0$
 \vspace{-0.3 cm}\[    
{\mathsf P}(Y(t)\le x|Y(t)>0)  = \int_0^\infty   {\mathsf P}(Z_{\lfloor{y}\rfloor}\le x|Z_{\lfloor{y}\rfloor}>0)\, d {\mathsf P}(N(t)\le y|Z_{N(t)}>0).
\]
\end{lemma}

{\bf Proof.} Using the independence of $Z_n$ and $N(t)$, we obtain
    \nbeq 
\lefteqn{{\mathsf P}(Y(t)\le x\ |\ Y(t)>0)  =  {\mathsf P}(Z_{N(t)} \le x\ |\ Z_{N(t)}>0)} \\
     & = & 
       \frac{1}{{\mathsf P}(Z_{N(t)}>0)} \sum_{n=0}^\infty {\mathsf P}(Z_{N(t)}\le x, Z_{N(t)}>0, N(t)=n) \nonumber \\
       & = & 
        \frac{1}{{\mathsf P}(Z_{N(t)}>0)} \sum_{n=0}^\infty {\mathsf P}(Z_{N(t)}\le x\ |\  Z_{N(t)}>0, N(t)=n) {\mathsf P}(Z_{N(t)}>0, N(t)=n) \nonumber \\ 
       & = & 
         \frac{1}{{\mathsf P}(Z_{N(t)}>0)} \sum_{n=0}^\infty {\mathsf P}(Z_n\le x\ |\  Z_n>0) {\mathsf P}(N(t)=n\ |\  Z_{N(t)}>0) {\mathsf P}(Z_{N(t)}>0) \nonumber \\
         & = & 
          \sum_{n=0}^\infty {\mathsf P}(Z_n\le x\ |\  Z_n>0) {\mathsf P}(N(t)=n\ |\  Z_{N(t)}>0) \nonumber \\
          & = & 
          \int_0^\infty {\mathsf P}(Z_{\lfloor{y}\rfloor}\le x\ |\ Z_{\lfloor{y}\rfloor}>0)\, d {\mathsf P}(N(t)\le y\ |\ Z_{N(t)}>0).\nonumber 
\neeq     

\hfill $\blacksquare$

{\bf Proof of Theorem 1(i)} 
It follows from (\ref{weights}) that for $\lfloor{x}\rfloor \le t <\lfloor{x+1}\rfloor$
\nbeq
{\mathsf P}(Y(t)>0) & = & 1- {\mathsf P}(Y(t)=0) \\
& = &  \sum_{k=1}^\infty \left[ {\mathsf P}(Z_k=0)-{\mathsf P}(Z_{k-1}=0)\right]{\mathsf P}(N(t)\le k)\\
& = &
\int_0^\infty {\mathsf P}(N(t)\le x)\, d (1-{\mathsf P}(Z{\lfloor{x}\rfloor}> 0)) \\
& = & - \int_0^\infty {\mathsf P}\left(  N(t)\le x\right)\, 
d {\mathsf P}(Z_{\lfloor{x}\rfloor}>0).
\neeq
Making the change of variables $y=x\mu/t$,  we obtain 
\[
\frac{\D {\mathsf P}(Y(t)>0)}{\D {\mathsf P}(Z_{\lfloor{t/\mu}\rfloor}>0)}
    =  
\int_0^\infty {\mathsf P}\left( \frac{N(t)}{t/\mu} \le y \right)\, 
d \frac{\D {\mathsf P}(Z_{\lfloor{yt/\mu}\rfloor}>0)}{\D {\mathsf P}(Z_{\lfloor{t/\mu}\rfloor}>0)}.
\]
Recall (\cite{BGT86}, Ch. 8.6.2) that if $\mu<\infty$, then  for $y\ge 1$ 
\be \label{nimu}
\lim_{t\to \infty}
{\mathsf P}\left(\frac{N(t)}{t/\mu}\le y\right)= 1.
\ee
Now, using Lemma \ref{Zlemma1} and (\ref{nimu}), we obtain
\nbeq
\lim_{t\to \infty} \frac{\D {\mathsf P}(Y(t)>0)}{\D {\mathsf P}(Z_{\lfloor{t/\mu}\rfloor}>0)}
    & = & 
- \lim_{t\to \infty}\int_0^\infty {\mathsf P}\left( \frac{N(t)}{t/\mu} \le y \right)\, 
d \frac{\D {\mathsf P}(Z_{\lfloor{yt/\mu}\rfloor}>0)}{\D {\mathsf P}(Z_{\lfloor{t/\mu}\rfloor}>0)}\\  
&  = & 
- \int_1^\infty \, d y^{-(1+|\theta|)}\\ 
    & = & 1.
\neeq
To complete the proof, it remains to observe that Lemma \ref{Zlemma1} implies 
\[{\mathsf P}(Z_{\lfloor{t/\mu}\rfloor}>0)\sim L_\theta(t)
\left(\frac{t}{\mu} \right)^{-(1+|\theta|)}, \quad t\to\infty.
\]

\hfill $\blacksquare$

{\bf Proof of Theorem 1(ii)}
Recalling the definition of the time to extinction $T$, we obtain  for $x\to \infty $
\[
W(x) = \sum_{k=0}^{\lfloor{x}\rfloor}w_k=\sum_{k=0}^{\lfloor{x}\rfloor} {\mathsf P}(T=k)={\mathsf P}(T\le \lfloor{x}\rfloor)\uparrow 1.
\]
Hence, by Lemma \ref{Zlemma2} with $0<\gamma<1$
\be \label{1-W}
1-W(x)= {\mathsf P}(T > \lfloor{x}\rfloor)={\mathsf P}(Z_{\lfloor{x}\rfloor}>0)\sim L_{\theta, \gamma}(x)x^{-\gamma}, \qquad x\to \infty.
\ee
Since
${\mathsf P}(Y(t)>0)=1-H(t)$, (\ref{lemma2i}), (\ref{1-W}) and Lemma \ref{l1} imply the claim.  

\hfill $\blacksquare$

{\bf Proof of Theorem 2} By Lemma \ref{integral}, making the change $u=y\mu/t$, we obtain
\beq \label{total_prob}
{\mathsf P}\left(\frac{Y(t)}{bt/\mu}\leq x\Big|Y(t)>0\right)\\
 & &  \hspace{-4cm}=  \int_{0}^{\infty }{\mathsf P}\left( \frac{Z_{\lfloor{y}\rfloor}}{ubt/\mu}\leq \frac{x}{u}
\Big|Z_{\lfloor{y}\rfloor}>0\right)\! \! d {\mathsf P}\! \left(N(t)\leq y\Big|Z_{N(t)}>0\right) \nonumber \\
&  & \hspace{-4 cm}= \int_{0}^{\infty }\! \! \! \! \! {\mathsf P}\left( \frac{Z_{\lfloor{t'}
\rfloor}}{bt'}\leq \frac{x}{u}
\Big|Z_{\lfloor{t'}\rfloor}>0\right) d {\mathsf P}\left(\frac{N(t)}{t/\mu}\leq u\Big| Z_{N(t)}>0\right), \nonumber
\eeq
where $t'=ut/\mu$.
According to Lemma \ref{Zlemma1} and Lemma \ref{Zlemma2}, we have for $x\ge 0$
\[ \label{lemma2ii22}
 \lim_{n\to \infty} \left(\frac{Z_n}{b n}\le x | Z_n>0\right)= {\mathsf P}(\delta \le x),
\]
where $\delta=\xi$ or $\delta=\eta_{\theta, \gamma}$, respectively.
Recall also (\ref{Dlimit})
\[  
\lim_{t\to\infty}{\mathsf P}\left( \frac{N(t)}{t/\mu}\le x | Z_{N(t)}>0\right)= \mathbf{1}_{\{x\ge 1 \}}.
\]
Therefore, taking the limit in (\ref{total_prob}) as $t\to \infty$ and applying the generalized Lebesgue dominated convergence theorem (see Theorem 2.4, \cite{S82}), we obtain
\[ 
\lim_{t\to \infty}{\mathsf P}\left(\frac{Y(t)}{bt/\mu}\leq x\Big|Y(t)>0\right)
= \int_{0}^\infty {\mathsf P}\left(\delta \le \frac{x}{u}\right) \, d\mathbf{1}_{\{u\ge 1\} } = {\mathsf P}(\delta\le x),
\]
which completes the proof.

\hfill $\blacksquare$

{\bf Proof of Theorem 3(i)} Under the assumptions, Lemma \ref{Zlemma1} implies that $\E[T]<\infty$. Therefore, it follows from Lemma \ref{Thm1_12}(i) and (\ref{Y and H}) that 
\be \label{pthm3i}
\lim_{t\to \infty} \frac{\D {\mathsf P}(Y(t)>0)}{\D {\mathsf P}(J>t)}=
\lim_{t\to \infty} \frac{\D 1-H(t)}{\D {\mathsf P}(J>t)}= \E[T-1].
\ee
Since ${\mathsf P}(J>t)\in {\cal RV}_{-\rho}$, relation (\ref{pthm3i}) completes the proof.  

\hfill $\blacksquare$

{\bf Proof of Theorem 3(ii)} Under the assumptions, Lemma \ref{Zlemma1} implies that $\E[T]=\infty$. It follows from Lemma \ref{Thm1_12}(ii) and (\ref{Y and H}) that 
\[
\lim_{t\to \infty} \frac{\D {\mathsf P}(Y(t)>0)}{\D {\mathsf P}(T>t^\rho/L_\rho(t))}
=\lim_{t\to \infty} \frac{\D 1-H(t)}{\D {\mathsf P}(T>t^\rho/L_\rho(t))}
=\frac{\Gamma^\gamma(2-\rho)\Gamma(1-\rho)}{\Gamma(1-\rho\gamma)}.
\]
To complete the proof, it remains to observe that (\ref{lemma2i}) implies that for $\gamma \in (0,1)$ we have
${\mathsf P}(T>n)={\mathsf P}(Z_n>0)\sim L_\gamma(n)n^{-\gamma}$ as $n\to \infty$.

\hfill $\blacksquare$

{\bf Proof of Theorem 4(i)} Similarly to (\ref{total_prob}), applying Lemma \ref{integral}, and making the change of variables $u=(y/t^\rho)L_\rho(t))$, we obtain
\[
{\mathsf P}\left(\frac{Y(t)}
{bt^\rho/L_\rho(t)}\leq x\Big|Y(t)>0\right) \! \! 
= \! \! \int_{0}^{\infty }\! \! \! \! \! {\mathsf P}\left( \frac{Z_{\lfloor{t''}
\rfloor}}{bt''}\leq \frac{x}{u}
\Big|Z_{\lfloor{t''}\rfloor}>0\right) d {\mathsf P}\left(\frac{N(t)}{t^\rho/L_\rho(t)}\leq u\Big| Z_{N(t)}>0\right), 
\]
where $t''=ut^\rho/L_\rho(t)$. Recall that according to Lemma \ref{Zlemma1}
\be \label{lemma2ii221}
 \lim_{n\to \infty} {\mathsf P}\left(\frac{Z_n}{b n}\le x | Z_n>0\right)= {\mathsf P}(\xi \le x), \qquad x\ge 0,
\ee
and also by  Lemma \ref{Thm2_12}(i)  
\be \label{lemma4i}
\lim_{t\to \infty}{\mathsf P}\left(\frac{N(t)}{t^\rho/L_\rho(t)}> x\ |\ T> N(t)\right) = \mathbf{1}_{\{x\ge 1 \}}.
\ee
Therefore, taking the limit  as $t\to \infty$ and using (\ref{lemma2ii221}) and (\ref{lemma4i}),  we obtain for $x\ge 0$
\[ 
\lim_{t\to \infty}{\mathsf P}\left(\frac{Y(t)}{bt/\mu}\leq x\Big|Y(t)>0\right)
= \int_{0}^\infty {\mathsf P}\left(\xi \le \frac{x}{u}\right) \, 
d\mathbf{1}_{\{u\ge 1 \}} = {\mathsf P}(\xi\le x).
\]

\hfill $\blacksquare$

{\bf Proof of Theorem 4(ii)} Similarly to (\ref{total_prob}), applying Lemma \ref{integral}, and making the change of variables $u=(y/t^\rho)L_\rho(t))$, we obtain
\[
{\mathsf P}\left(\frac{Y(t)}
{bt^\rho/L_\rho(t)}\leq x\Big|Y(t)>0\right) \! \! 
= \! \! \int_{0}^{\infty }\! \! \! \! \! {\mathsf P}\left( \frac{Z_{\lfloor{t''}
\rfloor}}{bt''}\leq \frac{x}{u}
\Big|Z_{\lfloor{t''}\rfloor}>0\right) d {\mathsf P}\left(\frac{N(t)}{t^\rho/L_\rho(t)}\leq u\Big| Z_{N(t)}>0\right), 
\]
where $t''=ut^\rho/L_\rho(t)$.  Recall that according to Lemma \ref{Zlemma2}
\[ \label{lemma2ii2}
 \lim_{n\to \infty} {\mathsf P}\left(\frac{Z_n}{b n}\le x | Z_n>0\right)= {\mathsf P}(\eta_{\theta, \gamma} \le x)\qquad x\ge 0,
\]
and also by Lemma \ref{Thm2_12}(ii) 
\[ \label{Dlimit3}
\lim_{t\to\infty}{\mathsf P}\left( \frac{N(t)}{t^\rho/L_\rho(t)}\le x | T> N(t)\right)= {\mathsf P}(\zeta_{\rho,\gamma}\le x) \qquad x\ge 0.
\]
Therefore,
\[
\hspace{-0.3cm}\lim_{t\rightarrow \infty }{\mathsf P}\left(\frac{ Y(t)}{bt^\rho/L_\rho(t)}\leq x\Big|Y(t)>0\right)
  =  \int_{0}^{\infty
}{\mathsf P}\left(\eta_{\theta, \gamma} \le \frac{x}{u}\right)\, d {\mathsf P}(\zeta_{\rho,\gamma}\le u)
   =   {\mathsf P}(\eta_{\theta, \gamma} \zeta_{\rho,\gamma}\le x).
\]
Here we used the fact that if two nonnegative r.v.s $X$ and $Y$ are independent and their d.f.s are $F_X$ and $F_Y$, respectively, then by the law of total probability for $z\ge 0$
\be \label{XY}
{\mathsf P}(XY\le z) =   {\mathsf P}\left(X\le \frac{z}{Y}\right) = \int_0^\infty 
F_X\left(\frac{z}{u}\right)\, dF_Y(u). 
\ee

\hfill $\blacksquare$

\section{Proofs of the Theorems for $\{U(t)\}_{t\ge 0}$}

The asymptotic behavior of $U(t)$ as $t\to \infty$ depends on the limiting behavior of its branching component $\{Y_k(t)\}_{t\ge 0}$. Assume that the following limits exist for $x\ge 0$
\be \label{Yk_limit}
\lim_{t\to \infty} \Pp \left(\frac{Y_k(t)}{M(t)}\le x\ | \ Y_k(t)>0\right) = F(x), \qquad k=1,2,\ldots,
\ee
where $M\in {\cal RV}_\rho$ for some $\rho\ge 0$.
Next lemma plays a key role in what follows.

\begin{lemma}[\cite{MY01}]
Suppose (\ref{Yk_limit}) holds. Assume 
 $\E[T]=\infty$ and
$
\Pp (Y(t)>0) \in {\cal RV_{-\beta}}$, $1/2 < \beta <1.
$

$\rm (i)$ If one of (\ref{alpha1}) or (\ref{alpha2}) is true and $0\le \Delta <\infty$, then for $x\ge 0$
\be \label{reg_lemmai}
\lim_{t\to \infty}\Pp\left(\frac{U(t)}{M(t)}\le x\right) = 
\left \{
\begin{array}{cl}
       \frac{\D \Delta}{\D \Delta +1} +\frac{\D G_1(x)}{\D \Delta +1} & \mbox{if} \ \ x >0 \\
       \frac{\D \Delta}{\D \Delta +1}     & \mbox{if} \ \ x=0,
            \end{array}
     \right. 
\ee
where $G_1$ is a proper d.f. given by
\be \label{G1}
G_1(x) = \frac{1}{{\rm B}(\beta, 1-\beta)} \int_0^1 F(xu^{-\rho})(1-u)^{\beta-1}u^{-\beta}\, du.
\ee
$\rm (ii)$ If (\ref{alpha2}) is true and $\Delta =\infty$, then for $x\ge 0$
\be \label{reg_lemmaii}
\lim_{t\to \infty}\Pp\left(\frac{U(t)}{M(t)}\le x \ |\ U(t)>0 \right) = G_2(x),
\ee
where $G_2$ is a proper d.f. given by
\be \label{G2}
G_2(x)=\frac{1}{{\rm B}(\alpha, 1-\beta)} \int_0^1 F(xu^{-\rho})(1-u)^{\alpha-1}u^{-\beta}\, du.
\ee
\end{lemma}
Let $\beta_{a,b}$ be a ${\rm Beta}(a,b)$ distributed r.v., independent of  
$\nu$, where $\nu$ has a d.f. $F$. Then, similarly to (\ref{XY}), using (\ref{G1}) and (\ref{G2}), we obtain that for $x\ge 0$
\be \label{XYkappa}
G_1(x) = \Pp(\nu \beta^\rho_{\beta, 1-\beta} \le x)\quad \mbox{and} \quad G_2(x) = 
\Pp(\nu \beta^\rho_{\alpha, 1-\beta} \le x). 
\ee

{\bf Proof of Theorem 5} Under the assumptions, Theorem 1(ii) implies  $\E [T]=\infty$. In Lemma 7 substitute
  $F(x)=\Pp(\eta_{\theta, \beta} \le x)$ and $\rho=1$.  Using (\ref{XYkappa}) with $\beta=\gamma$, we can rewrite (\ref{G1}) and (\ref{G2}) as
\be \label{G1_rewritten}
G_1(x) = \Pp(\eta_{\theta, \gamma} \beta_{\gamma, 1-\gamma} \le x)\quad \mbox{and} \quad G_2(x) = 
\Pp(\eta_{\theta, \gamma} \beta_{\alpha, 1-\gamma} \le x).
\ee
Here $\eta_{\theta, \gamma}$ is the limiting r.v. in Theorem 2(ii). On the other hand,
Theorem~1(ii) and Theorem~2(ii) yield 
$
{\mathsf P}(Y(t)>0)\in {\cal{RV}}_{-\gamma}$ for $1/2<\gamma<1 $,   
and  for $x\ge 0$ 
 \[
 \lim_{t\to \infty} {\mathsf P}\left(\frac{Y(t)}{b t/\mu}\le x | Y(t)>0\right)=
{\mathsf P}(\eta_{\theta, \gamma} \le x).
\]
Therefore, in view of (\ref{reg_lemmai}), (\ref{reg_lemmaii}), and (\ref{G1_rewritten}), the proof is complete. 

\hfill $\blacksquare$

{\bf Proof of Theorem 6} Under the assumptions,  Theorem 3(i) implies $\E [T]=\infty$. Substitute in Lemma 7,
  $F(x)=\Pp(\xi \le x)$ and $\beta=\rho$. Here $\xi$ is the limiting r.v. 
  in Theorem 4(i). Using (\ref{XYkappa}), we can rewrite (\ref{G1}) and (\ref{G2}) as
\be \label{G1_rewritten2}
G_1(x) = \Pp(\xi \beta^\rho_{\rho, 1-\rho} \le x)\quad \mbox{and} \quad G_2(x) =
\Pp(\xi \beta^\rho_{\alpha, 1-\rho} \le x).
\ee
Theorem~3(i) and Theorem~4(i) yield
$
{\mathsf P}(Y(t)>0)\in {\cal{RV}}_{-\rho}$ for $1/2<\rho<1
$
and thus
 \[
 \lim_{t\to \infty} {\mathsf P}\left(\frac{Y(t)}{t^\rho/L_\rho(t)}\le x | Y(t)>0\right)=
{\mathsf P}(\xi \le x).
\]
Therefore, referring to (\ref{reg_lemmai}), (\ref{reg_lemmaii}),  and (\ref{G1_rewritten2}), we complete the proof. 

\hfill $\blacksquare$

{\bf Proof of Theorem 7} Under the assumptions, Theorem 3(ii) implies $\E [T]=\infty$. In Lemma 7 substitute
  $F(x)=\Pp(\eta_{\theta, \gamma}\zeta_{\rho, \gamma} \le x)$ and $\beta=\rho\gamma$. Here $\eta_{\theta,\gamma}$
  and $\zeta_{\rho, \gamma}$ are the r.v.s 
  from Theorem 4(ii). Using (\ref{XYkappa}), we can rewrite (\ref{G1}) and (\ref{G2}) as
\be \label{G1_rewritten3}
G_1(x) = \Pp(\eta_{\theta, \gamma}\zeta_{\rho, \gamma} \beta^\rho_{\gamma, 1-\gamma} \le x)\quad \mbox{and} \quad G_2(x) =
\Pp(\eta_{\theta, \gamma}\zeta_{\rho, \gamma} \beta^\rho_{\alpha, 1-\gamma} \le x), \quad x\ge 0.
\ee
Recall that, $1/2<\rho \gamma<1$ and by
Theorems~3(ii) and 4(ii),
$
{\mathsf P}(Y(t)>0)\in {\mathcal{RV}}_{-\rho\gamma}$
and
 \[
 \lim_{t\to \infty} {\mathsf P}\left(\frac{Y(t)}{b t^\rho/L_\rho(t)}\le x | Y(t)>0\right)=
{\mathsf P}(\eta_{\theta,\gamma}\zeta_{\rho, \gamma} \le x), \qquad x\ge 0.
\]
Therefore, in view of (\ref{reg_lemmai}), (\ref{reg_lemmaii}),  and (\ref{G1_rewritten3}),  the proof is complete.

\hfill $\blacksquare$      

\section{Concluding Remarks}

The results in this study complement well the results presented in \cite{P20}. Here we  by focus on the regime of migration with dominating emigration ($\theta < 0$) instead of the case when the immigration dominates ($0 < \theta < 1$). 
The long-term limiting behavior of the subordinated branching process depends crucially on the tail 
  behavior of both interarrival time distribution and size of initial population. The normalizations in the limit theorems are affected by the interarrival time distribution, while the limiting 
  distribution itself retains memory of the size of the initial population (long memory effect of the initial condition). When the subordinator is a heavy-tailed renewal process, the classical theorems are replaced by results involving the time-changed processes and stable laws.

The paper also investigates the long-term behavior of the alternating regenerative branching process with migration, denoted by $\{U(t)\}_{t\ge 0}$. It alternates between active (branching) and inactive (sojourn at zero) periods. The main findings can be summarized as follows.
The asymptotic distribution of $U(t)$ as $t\to \infty$ depends on the tail behavior of both the interarrival time distribution (via the renewal process) and the sojourn time at zero, as well as the size of the initial population.
Consider the case when the mean interarrival time $\mu$ is finite and the initial population has infinite mean size. Then the normalized process $U(t)/(bt/\mu)$ converges in distribution to a mixture involving the product of  the subordinated process' limit and a Beta-distributed r.v., with an additional atom at zero.
Assume the mean interarrival time $\mu$ is infinite (heavy-tailed interarrival times). Then the normalization changes to $t^\rho/L_\rho(t)$, and the limiting distribution involves further randomization by a power-based Mittag-Leffler related r. v., reflecting the sublinear growth of the renewal process.
The parameter $\Delta$ (the asymptotic ratio of the tails of the sojourn time at zero and the length of the branching period) determines whether the limiting distribution has an atom at zero (if $0 \leq \Delta < \infty$) or is purely continuous (if $\Delta = \infty$).

Overall, the paper shows that the interplay 
between migration, heavy-tailed interarrival times, and 
regenerative structure leads to a variety of 
asymptotic behaviors, with explicit limiting distributions 
that generalize the classical results for the BGW processes.

It is worth exploring further the similarities and dissimilarities 
between Bellman-Harris processes and subordinated branching processes. Of particular interest would be a comparative study with the results for Bellman-Harris processes with infinite lifetime mean, large initial generation, and two-types particles (cf. \cite{Vatutin1979} and \cite{VatutinDyakonova2004}).

\section*{Acknowledgments}
The author gratefully acknowledges the support and constructive comments by Nikolay Yanev and Jordan Stoyanov, which helped improve paper's clarity and quality.

\bibliographystyle{plain}
\bibliography{main}

\end{document}